\documentclass[11pt]{article}
\usepackage{amscd}
\usepackage{amsfonts}
\usepackage{amsmath}
\usepackage{amssymb}
\usepackage{amsthm}
\usepackage{bbm}
\usepackage{CJK}
\usepackage{fancyhdr}
\usepackage{graphicx}
\usepackage{hyperref}
\usepackage{indentfirst}
\usepackage{latexsym}
\usepackage{mathrsfs}
\usepackage{xypic}

\usepackage[top=1in,bottom=1in,left=1.25in,right=1.25in]{geometry}
\newtheorem{theorem}{Theorem}[section]
\newtheorem{lemma}[theorem]{Lemma}
\newtheorem{definition}[theorem]{Definition}
\newtheorem{proposition}[theorem]{Proposition}
\newtheorem{example}[theorem]{Example}

\newtheorem{corollary}[theorem]{Corollary}

\usepackage[top=1in,bottom=1in,left=1.25in,right=1.25in]{geometry}
\def\<{\langle}
\def\>{\rangle}
\def\a{\alpha}
\def\b{\beta}

\def\c{\cdot}

\def\g{\gamma}

\def\o{\otimes}

\date{}
\begin{document}
\renewcommand{\baselinestretch}{1.2}
\renewcommand{\arraystretch}{1.0}
\title{\bf Representations and derivations of  Hom-Lie conformal superalgebras}
\author{{\bf Shuangjian Guo$^{1}$, Lihong Dong$^{2}$, Shengxiang Wang$^{3}$\footnote
        { Corresponding author(Shengxiang Wang):~~wangsx-math@163.com} }\\
{\small 1. School of Mathematics and Statistics, Guizhou University of Finance and Economics} \\
{\small  Guiyang  550025, P. R. of China} \\
{\small 2. College of Mathematics and Information Science, Henan  Normal University}\\
{\small  Xinxiang  453007, P. R. of  China}\\
{\small 3.~ School of Mathematics and Finance, Chuzhou University}\\
 {\small   Chuzhou 239000,  P. R. of China}}
 \maketitle
 \maketitle
\begin{center}
\begin{minipage}{13.cm}

{\bf \begin{center} ABSTRACT \end{center}}
 In this paper,  we introduce  a representation theory of Hom-Lie conformal superalgebras and  discuss the cases of adjoint representations. Furthermore,  we develop  cohomology theory of Hom-Lie conformal superalgebras and discuss
some applications to the study of deformations of regular Hom-Lie conformal superalgebras. Finally, we introduce derivations of multiplicative
Hom-Lie conformal superalgebras and study their properties.


 \end{minipage}
 \end{center}
 \normalsize\vskip1cm

\section{INTRODUCTION}
\def\theequation{0. \arabic{equation}}
\setcounter{equation} {0}

Lie conformal superalgebras are introduced by Kac in \cite{Kac98}, gives an axiomatic description of the singular part of the operator product expansion of chiral fields in conformal field theory. It is an useful tool to study vertex superalgebras and has many applications in the theory of Lie superalgebras. Moreover. Lie conformal superalgebras  have close connections to Hamiltonian formalism in the theory of nonlinear evolution equations.    Zhao, Yuan and Chen  developed deformation of Lie conformal superalgebras and   introduced derivations of multiplicative Lie conformal
superalgebras and study their properties in \cite{Zhao2017}.

In \cite{Ammar2010}, Ammar and Makhlouf introduced the notion of Hom-Lie superalgebras,
they gave a classification of Hom-Lie admissible superalgebras and proved a  graded version of Hartwig-Larsson-Silvestrov Theorem.
Later, Ammar, Makhlouf and Saadaoui \cite{Ammar2013} studied the representation and the cohomology of Hom-Lie superalgebras,
 and calculated the derivations and the second cohomology group of $q$-deformed Witt superalgebra.
In \cite{Yuan17}, Yuan introduced  Hom-Lie conformal superalgebra and proved that a Hom-Lie conformal superalgebra is equivalent to a Hom-Gel'fand-Dorfman superbialgebra.

  Recently, the Hom-Lie conformal algebra was introduced
and studied in \cite{Yuan14}, where it was proved that a Hom-Lie conformal algebra is equivalent to a Hom-Gel'fand-Dorfman bialgebra.   Zhao, Yuan and Chen  developed cohomology theory of Hom-Lie conformal algebras and discuss
some applications to the study of deformations of regular Hom-Lie conformal
algebras. Also, they  introduced  derivations of multiplicative Hom-Lie conformal
algebras and study their properties in \cite{Zhao2016}, which is different from \cite{Zhao2017}.

Motivated by these results, we introduce  a representation theory of Hom-Lie conformal superalgebras and  discuss the cases of adjoint representations. Furthermore,  we develop  cohomology theory of Hom-Lie conformal superalgebras and discuss
some applications to the study of deformations of regular Hom-Lie conformal superalgebras. Finally, we introduce derivations of multiplicative
Hom-Lie conformal superalgebras and study their properties.

This paper is organized as follows.

In Section 2,  we introduce  a representation theory of Hom-Lie conformal superalgebras and  discuss the cases of adjoint representations.

In Section 3, we develop  cohomology theory of Hom-Lie conformal superalgebras and discuss
some applications to the study of deformations of regular Hom-Lie conformal superalgebras.

In Section 4, we study derivations of multiplicative Hom-Lie conformal superalgebras.
and prove  the direct sum of the space of derivations is a Hom-Lie
conformal superalgebra. In particular, any derivation gives rise to a derivation extension of
a multiplicative Hom-Lie conformal superalgebra.

In Section 5, we introduce  generalized derivations of multiplicative
Hom-Lie conformal superalgebras and study their properties.

Throughout the paper, all algebraic systems are supposed to be over a field $\mathbb{C}$, and denote by $\mathbb{Z}_{+}$ the set of all nonnegative integers and by $\mathbb{Z}$.

\section{Representations of Hom-Lie conformal superalgebras }
\def\theequation{\arabic{section}. \arabic{equation}}
\setcounter{equation} {0}

In this section, we introduce  a representation theory of Hom-Lie conformal superalgebras and  discuss the cases of adjoint representations.

\begin{definition}$^{\cite{Yuan17}}$
A Hom-Lie conformal superalgebra $R=R_{\bar{0}}\oplus R_{\bar{1}}$ is a  $\mathbb{Z}_2$-graded $\mathbb{C}[\partial]$-module equipped with an even linear endomorphism
$\a$ such that $\a\partial=\partial\a$, and  a $\mathbb{C}$-linear map
\begin{eqnarray*}
R\o R\rightarrow \mathbb{C}[\lambda]\o R,  ~~~a\o b\mapsto [a_\lambda b]
\end{eqnarray*}
such that $[R_{\phi\lambda}R_{\varphi}]\subseteq R_{\phi+\varphi}[\lambda]$, $\phi,\varphi\in \mathbb{Z}_2$, and the following axioms hold for $a,b,c\in R$
\begin{eqnarray}
&&[\partial a_{\lambda}b] =-\lambda[a_\lambda b],[ a_{\lambda}\partial b] =(\partial+\lambda)[a_\lambda b],\\
&&[a_{\lambda}b]=-(-1)^{|a|||b|}[b_{-\lambda-\partial}a],\\
&&[\a(a)_\lambda[b_\mu c]]=[[a_{\lambda}b]_{\lambda+\mu}\a(c)]+(-1)^{|a|||b|}[\a(b)_{\mu}[a_{\lambda}c]].
\end{eqnarray}
\end{definition}
A Hom-Lie conformal superalgebra $(R, \a)$ is called finite if $R$ is a finitely generated $\mathbb{C}[\partial]$-module.
The rank of $R$ is its rank as a $\mathbb{C}[\partial]$-module.

A Hom-Lie conformal superalgebra $(R, \a)$  is called multiplicative if $\a$ is an algebra endomorphism, i.e., $\a([a_\lambda b])=[\a(a)_{\lambda}\a(b)]$ for any $a,b\in R$. In particular, if $\a$ is an algebra isomorphism, then $(R, \a)$  is regular.

\begin{example}
Let $g=g_{\bar{0}}\oplus g_{\bar{1}}$ be a complex Hom-Lie superalgebra with Lie bracket $[-,-]$. Let $(Curg)_{\theta}:=\mathbb{C}[\partial]\o g_\theta$ be the free $\mathbb{C}[\partial]$-module. Then $Curg=(Curg)_{\bar{0}}\oplus (Curg)_{\bar{1}}$ is a Hom-Lie conformal superalgebra with $\lambda$-bracket given by
\begin{eqnarray*}
&& \a(f(\partial)\o a)=f(\partial)\o \a(a),\\
&&[(f(\partial)\o a)_{\lambda}(g(\partial)\o b)]=f(-\lambda)g(\partial+\lambda)\o [a,b], \forall a,b\in g.
\end{eqnarray*}
\end{example}

\begin{example}
Let $R=\mathbb{C}[\partial]L\oplus \mathbb{C}[\partial]E$ be a free $\mathbb{Z}_2$-graded $\mathbb{C}[\partial]$-module. Define
\begin{eqnarray*}
&& \a(L)=f(\partial)L, \a(E)=g(\partial)E,\\
&& [L_\lambda L]=(\partial+2\lambda)L, [L_\lambda E]=(\partial+\frac{3}{2}\lambda)E, [E_\lambda L]=(\frac{1}{2}\partial+\frac{3}{2}\lambda)E, [E_\lambda E]=0,
\end{eqnarray*}
where $R_{\bar{0}}=\mathbb{C}[\partial]L$ and $ R_{\bar{1}}= \mathbb{C}[\partial]E$. Then $(R,\a)$ is a Hom-Lie conformal superalgebra.
\end{example}

\begin{definition}
Let $(M, \b)$ and $(N, \omega)$ be $\mathbb{Z}_2$-graded $\mathbb{C}[\partial]$-modules. A Hom-conformal linear map of degree $\theta$ from $M$ to $N$ is a sequence $f={f_{(n)}}_{n\in \mathbb{Z}_{\geq 0}}$ of $f_{(n)}\in Hom_{\mathbb{C}}(M, N)$ satisfying that
\begin{eqnarray*}
&&\partial_N f_{(n)}-f_{(n)}\partial_M=-n f_{(n-1)}, f_{\lambda}(M_\mu)\subseteq N_{\mu+\theta}, n\in \mathbb{Z}_{\geq 0}, \mu, \theta\in \mathbb{Z}_2,\\
&& \partial_M\b=\b\partial_M,  \partial_N\omega=\omega\partial_N, f_{(n)}\b=\omega f_{(n)}.
\end{eqnarray*}
Set $f_\lambda=\sum_{n=0}^{\infty}\frac{\lambda^{n}}{n!}f_{(n)}$. Then $f$ is a Hom-conformal linear map of degree $\theta$ if and only if
\begin{eqnarray*}
&&f_{\lambda} \partial_M=(\partial_N+\lambda)f_\lambda, f_{\lambda}(M_\mu)\subseteq N_{\mu+\theta}[\lambda],\\
&& \partial_M\b=\b\partial_M,  \partial_N\omega=\omega\partial_N, f_{\lambda}\b=\omega f_{\lambda}.
\end{eqnarray*}
\end{definition}

\begin{definition}
A Hom-associative conformal superalgebra $R$ is a  $\mathbb{Z}_2$-graded $\mathbb{C}[\partial]$-module equipped with an even linear endomorphism
$\a$ such that $\a\partial=\partial\a$, endowed with a $\lambda$-product from $R\o R$ to $\mathbb{C}[\partial]\o R$, for any $a,b,c\in R$, satisfying the following conditions:
\begin{eqnarray}
&&(\partial a)_\lambda b=-\lambda a_{\lambda}b, a_{\lambda}(\partial b)=(\partial+\lambda)a_{\lambda}b,\nonumber\\
&& \a(a)_{\lambda}(b_{\mu}c)=(a_\lambda b)_{\lambda+\mu}\a(c).
\end{eqnarray}
\end{definition}

\begin{theorem}
Let $(R,\alpha)$ be a Hom-associative conformal superalgebra with an even linear endomorphism $\alpha$.
One can define
 \begin{eqnarray*}
[a_{\lambda}b]=a_{\lambda} b-(-1)^{|a||b|}b_{-\lambda-\partial} a, \forall a,b\in R.
  \end{eqnarray*}
 Then $(R, \alpha)$ is a  Hom-Lie conformal  superalgebra.

{\bf Proof.} We only prove (2.3). For any $a,b,c\in R$, we have
\begin{eqnarray*}
&&[\a(a)_{\lambda}[b_{\mu}c]]\\
&=&[\a(a)_{\lambda}(b_{\mu}c-(-1)^{|b||c|}c_{-\mu-\partial} b)]\\
&=&\a(a)_{\lambda}(b_{\mu}c)-(-1)^{|a|(|b|+|c|)}(b_{\mu}c)_{-\lambda-\partial}\a(a)\\
&&-(-1)^{|b||c|}\a(a)_{\lambda}(c_{-\mu-\partial} b)+(-1)^{|a||c|+|a||b|+|b||c|}(c_{-\mu-\partial} b)_{-\lambda-\partial}\a(a).
  \end{eqnarray*}
Similarly, we have
\begin{eqnarray*}
&&[[a_\lambda b]_{\lambda+\mu}\a(c)]\\
&=& [(a_{\lambda} b-(-1)^{|a||b|}b_{-\lambda-\partial} a)_{_{\lambda+\mu}}\a(c)]\\
&=& (a_{\lambda} b)_{\lambda+\mu}\a(c)-(-1)^{|c|(|a|+|b|)}\a(c)_{-\lambda-\mu-\partial}(a_{\lambda} b)\\
&&-(-1)^{|a||b|}(b_{-\lambda-\partial} a)_{_{\lambda+\mu}}\a(c)+(-1)^{|a||c|+|a||b|+|b||c|}\a(c)_{-\lambda-\mu-\partial}(b_{-\lambda-\partial} a)_{\lambda+\mu}\\
&&(-1)^{|a||b|}[\a(b)_{\mu},[a_{\lambda}c]]\\
&=&(-1)^{|a||b|}\a(b)_{\mu}(a_{\lambda}c)-(-1)^{|a|(|c|+|b|)}\a(b)_{\mu}(c_{-\lambda-\partial} a)\\
&&-(-1)^{|a||c|}(a_{\lambda}c)_{-\mu-\partial}\a(b)+(-1)^{|a||b|}(c_{-\lambda-\partial} a)_{-\mu-\partial}\a(b).
  \end{eqnarray*}
By the associativity  (2.4), it is not hard to check that
\begin{eqnarray*}
[\a(a)_\lambda[b_\mu c]]=[[a_{\lambda}b]_{\lambda+\mu}\a(c)]+(-1)^{|a|||b|}[\a(b)_{\mu}[a_{\lambda}c]],
\end{eqnarray*}

as desired. And this finishes the proof.  \hfill $\square$

\end{theorem}

Let Chom$(M, N)_{\theta}$ denote the set of Hom-conformal linear maps of degree $\theta$ from $M$ to $N$. Then  Chom$(M, N)$= Chom$(M, N)_{\bar{0}}\oplus$  Chom$(M, N)_{\bar{1}}$ is a $\mathbb{Z}_2$-graded $\mathbb{C}[\partial]$-module via
\begin{eqnarray*}
\partial f_{(n)}=-nf_{(n-1)},  ~~~\mbox{equivalently}~~~  \partial f_{\lambda}=-\lambda f_{\lambda}.
\end{eqnarray*}
The composition $f_\lambda g:L\rightarrow N\o \mathbb{C}[\lambda]$ of Hom-conformal linear maps $f: M\rightarrow N$ and $g: L\rightarrow M$ is given by
\begin{eqnarray*}
(f_\lambda g)_{\lambda+\mu}=f_\lambda g_\mu, \forall f,g\in  Chom(M, N).
\end{eqnarray*}
If $(M, \b)$ is a finitely generated $\mathbb{Z}_2$-graded $\mathbb{C}[\partial]$-module, then Cend$(M):=$Chom$(M, N)$ is a Hom-associative
conformal superalgebra with respect to the above composition. Thus, Cend$(M)$ becomes a
Hom-Lie conformal superalgebra, denoted as $gc(M)$, with respect to the following $\lambda$-bracket
\begin{eqnarray*}
[f_\lambda g] =f_\lambda g_{\mu-\lambda}-(-1)^{|f||g|}g_{\mu-\lambda}f_{\lambda}.
\end{eqnarray*}
\begin{definition}
Let $(R,\a)$ be a Hom-Lie conformal superalgebra. A representation of  $R$ is a triple $(\rho, M, \b)$, where $(M, \b)$ is a  $\mathbb{Z}_2$-graded $\mathbb{C}[\partial]$-module,  and   $\rho: R\rightarrow Cend{M}$ is a  linear map  satisfying
\begin{eqnarray}
&&\rho(\partial(a))_{\lambda}=-\lambda\rho(a)_{\lambda},\rho\partial=\partial\rho, [\rho(a)_{\lambda}, \rho(b)_{\mu}]=\rho([a_\lambda b])_{\lambda+\mu},\nonumber\\
&&\rho([a_{\lambda}b])_{\lambda+\mu}\b=\rho(\a(a))_{\lambda}\rho(b)_{\mu}-(-1)^{|a||b|}\rho(\a(b))_{\mu}\rho(a)_{\lambda}.
\end{eqnarray}

\end{definition}

\begin{proposition}
Let $(R,\a)$ be a Hom-Lie conformal superalgebra  and $(\rho, M, \b)$ is a representation of the Hom-Lie conformal superalgebra $R$. The direct sum $R\oplus M$ with a $\lambda$-bracket $[\c_{\lambda}\c]$ on $R\oplus M$ defined by
\begin{eqnarray*}
[(a+u)_{\lambda}(b+v)]_{M}=[a_{\lambda}b]+\rho(a)_{\lambda}v-(-1)^{|a||b|}\rho(b)_{-\partial-\lambda}u,~~~\forall a,b\in R, u,v\in M.
\end{eqnarray*}
and the twist map $\a+\b: R\oplus M\rightarrow R\oplus M$ defined by
\begin{eqnarray*}
(\a+\b)(a+u)=\a(a)+\b(u),  \forall a\in R, u\in M
\end{eqnarray*}
is a Hom-Lie conformal superalgebra.
\end{proposition}
{\bf Proof.} Note that $R\oplus M$ is equipped with a $\mathbb{C}[\partial]$-module structure via
\begin{eqnarray*}
\partial(a+u)=\partial(a)+\partial(u), \forall a\in R, u\in M.
\end{eqnarray*}
It is easy to check that $(\a+\b)\circ \partial=\partial\circ (\a+\b)$, Now, we check that (2.1) holds, for any $a,b\in R$ and $u,v\in M$, we have
\begin{eqnarray*}
[\partial(a+u)_{\lambda}(b+v)]_{M}&=&[(\partial a+\partial u)_{\lambda}(b+v)]_{M}\\
&=& [\partial a_{\lambda}b]+\rho(\partial a)_{\lambda}v-(-1)^{|a||b|}\rho(b)_{-\partial-\lambda}\partial u\\
&=& -\lambda[a_{\lambda}b]-\lambda\rho(a)_{\lambda}v-(-1)^{|a||b|}(\partial-\lambda-\partial)\rho(b)_{-\partial-\lambda}u\\
&=&  -\lambda([a_{\lambda}b]+\rho(a)_{\lambda}v-(-1)^{|a||b|}\rho(b)_{-\partial-\lambda}u)\\
&=&  -\lambda [(a+u)_{\lambda}(b+v)]_{M},
\end{eqnarray*}
and
\begin{eqnarray*}
[(a+u)_{\lambda}\partial(b+v)]_{M}&=& [(a+u)_{\lambda}(\partial b+\partial v)]_{M}\\
&=& [a_\lambda \partial b]+\rho( a)_{\lambda}\partial v-(-1)^{|a||b|}\rho(\partial b)_{-\partial-\lambda}u\\
&=&(\partial+\lambda) [a_\lambda  b]+(\partial+\lambda)\rho( a)_{\lambda} v-(-1)^{|a||b|}(\partial+\lambda)\rho(b)_{-\partial-\lambda}u\\
&=& (\partial+\lambda)([a_{\lambda}b]+\rho(a)_{\lambda}v-(-1)^{|a||b|}\rho(b)_{-\partial-\lambda}u)\\
&=& (\partial+\lambda)[(a+u)_{\lambda}(b+v)]_{M}.
\end{eqnarray*}
(2.2) follows from
\begin{eqnarray*}
[(b+v)_{-\partial-\lambda}(a+u)]_{M}&=&[b_{-\partial-\lambda}a]+\rho(b)_{-\partial-\lambda}u-(-1)^{|a||b|}\rho(a)_{\lambda}v\\
&=&-(-1)^{|a||b|}([a_{\lambda}b]-(-1)^{|a||b|}\rho(b)_{-\partial-\lambda}u+\rho(a)_{\lambda}v)\\
&=& -(-1)^{|a||b|}[(a+u)_{\lambda}(b+v)]_{M}.
\end{eqnarray*}
Next we show that the Hom-Jacobi identity, we compute
\begin{eqnarray*}
&&[(\a+\b)(a+u)_{\lambda}[(b+v)_{\mu}(c+w)]_M]_M\nonumber\\
&=& [(\a(a)+\b(u))_{\lambda}[(b+v)_{\mu}(c+w)]_M]_M\nonumber\\
&=& [(\a(a)+\b(u))_{\lambda}[[b_{\mu}c]+\rho(b)_{\mu}w-(-1)^{|b||c|}\rho(c)_{-\partial-\mu}v]_M\nonumber\\
&=& [\a(a)_{\lambda}[b_{\mu}c]]+\rho(\a(a))_{\lambda}(\rho(b)_{\mu}w)-(-1)^{|c||b|}\rho(\a(a))_{\lambda}(\rho(c)_{-\partial-\mu}v)\nonumber\\
&&-(-1)^{(|c|+|b|)|a|}\rho([b_\mu c])_{-\partial-\lambda}(\b(u)),~~~~~~
\end{eqnarray*}
\begin{eqnarray*}
&&(-1)^{|a||b|}[(\a+\b)(b+v)_{\mu}[(a+u)_{\lambda}(c+w)]_M]_M\nonumber\\
&=& (-1)^{|a||b|}[\a(b)_{\mu}[a_{\lambda}c]]+(-1)^{|a||b|}\rho(\a(b))_{\mu}(\rho(a)_{\lambda}w)\nonumber\\
&&-(-1)^{|a||b|+|c||a|}\rho(\a(b))_{\mu}(\rho(c)_{-\partial-\lambda}u)-(-1)^{|c||b|}\rho([a_{\lambda}c]_{-\partial-\mu})\b(v)~~~~~~
\end{eqnarray*}
and
\begin{eqnarray*}
&&[[(a+u)_{\lambda}(b+v)]_{M_{\lambda+\mu}}(\a+\b)(c+w)]_M\nonumber\\
&=& [[a_{\lambda}b]_{\lambda+\mu}\a(c)]+ \rho([a_{\lambda}b])_{\lambda+\mu}\b(w)-(-1)^{|c|(|a|+|b|)}\rho(\a(c))_{-\partial-\lambda-\mu}(\rho(a)_{\lambda}v)\nonumber\\ &&+(-1)^{|c|(|b|+|a|)+|b||a|}\rho(\a(c))_{-\partial-\lambda-\mu}(\rho(b)_{-\partial-\lambda}u)~~~~~~~
\end{eqnarray*}

Since $(\rho, M, \b)$ is a representation of the Hom-Lie conformal superalgebra $R$, we have
\begin{eqnarray*}
\rho(\a(a))_{\lambda}(\rho(b)_{\mu}w)-(-1)^{|a||b|}\rho(\a(b))_{\mu}(\rho(a)_{\lambda}w)=\rho([a_{\lambda}b])_{\lambda+\mu}\b(w).
\end{eqnarray*}
Thus  $(R\oplus M, \a+\b)$ is a Hom-Lie conformal superalgebra. \hfill $\square$

\begin{proposition}
Let $(R,\a)$ be a Hom-Lie conformal superalgebra  and a Hom-conformal linear map $ad: R\rightarrow  Cend{R}$  defined by $(ad (a))_{\lambda}(b)=[a_{\lambda}b]$.
Then $(R, ad, \a)$ is a representation of $R$.
\end{proposition}
{\bf Proof.} Since $(R,\a)$ is a Hom-Lie conformal superalgebra, the Hom-Jacobi identity (2.3) may be written for any $a,b,c\in R$
\begin{eqnarray*}
ad([a_\lambda b])_{\lambda+\mu}(\a(c))=ad(\a(a))_{\lambda}(ad(b)_{\mu}(c))-(-1)^{|a||b|}ad(\a(b))_{\mu}(ad(a)_{\lambda}(c))
\end{eqnarray*}
Then the conformal linear map $ad$ satisfies
\begin{eqnarray*}
ad([a_\lambda b])_{\lambda+\mu} \a=ad(\a(a))_{\lambda} ad(b)_{\mu}-(-1)^{|a||b|}ad(\a(b))_{\mu} ad(a)_{\lambda}
\end{eqnarray*}
Therefore, $(R, ad, \a)$ is a representation of $R$.   \hfill $\square$

We call the representation defined in the previous proposition adjoint representation of the Hom-Lie conformal superalgebra.

In the following, we explore the dual representation and coadjoint representation of  Hom-Lie conformal superalgebras.

Let $(R,\a)$ be a Hom-Lie conformal superalgebra and $(\rho, M, \b)$ be a representation of $R$. Let $M^{\ast}$ be the dual vector space of $M$. We define a Hom-conformal linear map $\tilde{\rho}: R\rightarrow Cend{M^{\ast}}$ by $\tilde{\rho}_{\lambda}(a)=-\rho_{\lambda}(x)$.

Let $f\in M^{\ast}, a,b\in R$ and $u\in M$. We compute the right hand side of (2.4)
\begin{eqnarray*}
&&(\tilde{\rho}(\a(a))_{\lambda})\tilde{\rho}(b)_{\mu}-(-1)^{|a||b|}(\tilde{\rho}(\a(b))_{\mu})\tilde{\rho}(a)_{\lambda} (f)(u)\\
&=& (\tilde{\rho}(\a(a))_{\lambda})(\tilde{\rho}(b)_{\mu}(f))-(-1)^{|a||b|}(\tilde{\rho}(\a(b))_{\mu})(\tilde{\rho}(a)_{\lambda}(f)) (u)\\
&=&  -(\tilde{\rho}(b)_{\mu}(f))(\rho(\a(a))_{\lambda}(u))+(-1)^{|a||b|}(\tilde{\rho}(a)_{\lambda}(f))(\tilde{\rho}(\a(b))_{\mu} (u))\\
&=& f(\rho(b)_{\mu}\rho(\a(a))_{\lambda}(u))-(-1)^{|a||b|}f(\rho(a)_{\lambda}\rho(\a(b))_{\mu}(u))\\
&=& f((\rho(b)_{\mu}\rho(\a(a))_{\lambda}-(-1)^{|a||b|}f(\rho(a)_{\lambda}\rho(\a(b))_{\mu})(u)).
\end{eqnarray*}
On the other hand, we set that the map $\tilde{\b}=\b$, then the left  hand side of (2.4),
\begin{eqnarray*}
((\tilde{\rho}([a_\lambda b])_{\lambda+\mu}\tilde{\b}) (f))(u)=((\tilde{\rho}([a_\lambda b])_{\lambda+\mu} (f\b))(u))=-f \b(\rho([a_\lambda b])_{\lambda+\mu}(u)).
\end{eqnarray*}
Therefore, we have the following proposition:
\begin{proposition}
Let $(R,\a)$ be a Hom-Lie conformal superalgebra and $(\rho, M, \b)$ be a representation of $R$.
the $(V^{\ast},\tilde{\rho}, \tilde{\b} )$, where  $\tilde{\rho}: R\rightarrow Cend{M^{\ast}}$ is given by $\tilde{\rho}_{\lambda}(a)=-\rho_{\lambda}(x)$, defines a representation of  Hom-Lie conformal superalgebra $(R,\a)$ if and only if
\begin{eqnarray*}
\b\rho([a_\lambda b])_{\lambda+\mu}=(-1)^{|a||b|}\rho(a)_{\lambda}\rho(\a(b)_{\mu}-   \rho(b)_{\mu}\rho(\a(a))_{\lambda}, \forall a,b\in R.
\end{eqnarray*}
\end{proposition}

\begin{corollary}
Let $(R,\a)$ be a Hom-Lie conformal superalgebra and $(R, ad, \a)$ be the adjoint representation of $R$, where $ad: g\rightarrow Cend{g}$. We set $\tilde{ad}: g\rightarrow Cend{g^{\ast}}$ and $\tilde{ad}(x)_{\lambda}(f)=-f ad(x)_{\lambda}$. Then $(g^{\ast},\tilde{ad}, \tilde{\a} )$  is a representation of $R$ if and only if
\begin{eqnarray*}
\a([a_{\lambda}b]_{\lambda+\mu}, c)=(-1)^{|a||b|}[a_{\lambda} [\a(b)_{\mu}c]]-[y_{\mu}[a_\lambda c]], \forall a,b,c\in R.
\end{eqnarray*}
\end{corollary}

\section{ Cohomology  of Hom-Lie conformal superalgebras }
\def\theequation{\arabic{section}. \arabic{equation}}
\setcounter{equation} {0}

In this section, we introduce the cohomology theory of Hom-Lie conformal superalgebras
and  the Nijenhuis operators of Hom-Lie conformal superalgebras. As an application of the
cohomology theory, we give the definition of a deformation and show that the deformation generated
by a 2-cocycle Nijenhuis operator is trivial.

\begin{definition}
An $n$-cochain ($n\in \mathbb{Z}_{\geq0}$) of a regular Hom-Lie conformal superalgebra $R$ with with coefficients in a module $(M,\b )$ is a $\mathbb{C}$-linear map
\begin{eqnarray*}
&&\g:R^n\rightarrow M[\lambda_1,... ,\lambda_n] \\
&&(a_1,...,a_n)\mapsto \g_{\lambda_1,...,\lambda_n}(a_1,... , a_n).
\end{eqnarray*}
where $M[\lambda_1,...,\lambda_n]$ denotes the space of polynomials with coefficients in $M$, satisfying the following conditions:

(1) Conformal antilinearity:
\begin{eqnarray*}
\g_{\lambda_1,...,\lambda_n}(a_1,... ,\partial a_i,... ,a_n)=-\lambda_i\g_{\lambda_1,...,\lambda_n}(a_1,... ,a_i,... ,a_n).
\end{eqnarray*}

(2) Skew-symmetry:
\begin{eqnarray*}
&&\gamma(a_{1},..., a_{i}, a_{j},..., a_{n})=-(-1)^{|a_i||a_j|}\gamma(a_{1}, ..., a_{j}, a_{i}, ... , a_{n}).
\end{eqnarray*}

(3)Commutativity:
\begin{eqnarray*}
\g\circ \a=\b\circ \g
\end{eqnarray*}
\end{definition}

Let $R^{\otimes 0} = \mathbb{C}$ as usual so that a $0$-cochain
 is an element of $M$. Define a differential $d$ of a cochain $\g$ by
\begin{eqnarray*}
&&(d\gamma)_{\lambda_1,...,\lambda_{n+1}}(a_{1},...,a_{n+1})\nonumber\\
&=&\sum_{i=1}^{ n+1}(-1)^{i+1}(-1)^{(|\g|+|a_1|+...|a_{i-1}|)|a_i|}\rho(\alpha^{n}(a_{i}))_{\lambda_i}\gamma_{\lambda_1,...,\hat{\lambda}_{i},...,\lambda_{n+1}}(a_{1},..., \hat{a}_{i},..., a_{n+1})\\
  &&+\sum_{1\leq i<j\leq n+1} (-1)^{i+j}(-1)^{(|a_1|+...|a_{i-1}|)|a_i|+(|a_i|+...+|a_{j-1}|)||a_j+|a_i||a_j|}\\
  &&\gamma_{\lambda_i+\lambda_j, \lambda_1,...,\hat{\lambda}_i,...\hat{\lambda}_j,...,\lambda_{n+1}}([a_{i\lambda_i}a_{j}], \a(a_{1}),... ,\hat{a}_{i},... , \hat{a}_{j},...,\a(a_{n+1}))
\end{eqnarray*}
where $\g$ is extended linearly over the polynomials in
$\lambda_i$. In particular, if $\g$ is a 0-cochain, then $(d\g)_\lambda a=a_\lambda \g$.

\begin{proposition}
$d\g$ is a cochain and $d^2=0$.
\end{proposition}

 {\bf Proof.} Let $\g$ be an $n$-cochain.  It is easy to check that $d$ satisfies conformal antilinearity and skew-symmetry. Commutativity is obviously satisfied. Thus $d$ is an $(n + 1)$-cochain.

 Next we will check that $d^2=0$.  By straightforward computation, we have
 \begin{eqnarray}
&& (d^2\g) _{\lambda_1,...,\lambda_{n+2}}(a_{1},...,a_{n+2})\nonumber\\
&=& \sum_{i=1}^{ n+1}(-1)^{i+1+|\g||a_i|+|A_i|}\alpha^{n+1}(a_{i})_{\lambda_i}(d\g)_{\lambda_1,...,\hat{\lambda}_{i},...,\lambda_{n+2}}(a_{1},..., \hat{a}_{i},..., a_{n+2})\nonumber\\
 && +\sum_{1\leq i<j\leq n+1} (-1)^{i+j+A_i+A_j+|a_i||a_j|}(d\g)_{\lambda_i+\lambda_j, \lambda_1,...,\hat{\lambda}_i,...\hat{\lambda}_j,...,\lambda_{n+2}}\nonumber\\
  &&([\a^{-1} \a(a_{i})_{\lambda_i}a_{j}], \a(a_{1}),... ,\hat{a}_{i},... , \hat{a}_{j},...,\a(a_{n+2}))\nonumber\\
  &=&\sum_{i=1}^{n+2}\sum_{j=1}^{i-1}(-1)^{i+j+A_i+A_j+|\g|(|a_i|+|a_j|)}\nonumber\\
  &&\rho(\alpha^{n+1}(a_{i}))_{\lambda_i}(\rho(\alpha{n+1}(a_{j}))_{\lambda_j}\g_{\lambda_1,...,\hat{\lambda}_{j,i},...,\lambda_{n+2}}(a_{1},..., \hat{a}_{j,i},..., a_{n+2})\\
  &&\sum_{i=1}^{n+2}\sum_{j=i+1}^{n+2}(-1)^{i+j+1+A_i+(A_j-|a_i|)+|\g|(|a_i|+|a_j|}\nonumber\\
  &&\rho(\alpha^{n+1}(a_{i}))_{\lambda_i}(\rho(\alpha^{n}(a_{j}))_{\lambda_j}\g_{\lambda_1,...,\hat{\lambda}_{i,j},...,\lambda_{n+2}}(a_{1},..., \hat{a}_{i,j},..., a_{n+2})\\
  &&+\sum_{i=1}^{n+2}\sum_{1\leq j< k<i}^{n+2}(-1)^{i+j+k+1+|\g||a_i|+A_i+A_j+A_k+|a_j||a_k|}\rho(\alpha^{n+1}(a_{i}))_{\lambda_i}\nonumber\\
  &&\g_{\lambda_j+\lambda_k, \lambda_1,...,\hat{\lambda}_{j,k,i},...,\lambda_{n+2}}([a_{j\lambda_j}a_{k}], \a(a_{1}),... ,\hat{a}_{j,k,i},...,\a(a_{n+2}))\\
  &&+\sum_{i=1}^{n+2}\sum_{1\leq j< i<k}^{n+2}(-1)^{i+j+k+1+|\g||a_i|+A_i+A_j+(A_k-|a_i|)+|a_j||a_k|}\rho(\alpha^{n+1}(a_{i}))_{\lambda_i}\nonumber\\
  &&\g_{\lambda_j+\lambda_k, \lambda_1,...,\hat{\lambda}_{j,i,k},...,\lambda_{n+2}}([a_{j\lambda_j}a_{k}], \a(a_{1}),... ,\hat{a}_{j,i,k},...,\a(a_{n+2}))
   \end{eqnarray}
   \begin{eqnarray}
   &&+\sum_{i=1}^{n+2}\sum_{1\leq i< j<k}^{n+2}(-1)^{i+j+k+1+|\g||a_i|+A_i+A_j+(A_j-|a_i|)+(A_i-|a_i|)+|a_j||a_k|}\rho(\alpha^{n+1}(a_{i}))_{\lambda_i}\nonumber\\
  &&\g_{\lambda_j+\lambda_k, \lambda_1,...,\hat{\lambda}_{i,j,k},...,\lambda_{n+2}}([a_{j\lambda_j}a_{k}], \a(a_{1}),... ,\hat{a}_{i,j,k},...,\a(a_{n+2}))\\
  &&+\sum_{1\leq i< j}^{n+2}\sum_{k=1}^{i-1}(-1)^{i+j+k+A_i+A_j+A_k+|a_i||a_j|+(|\g|+|a_i|+|a_j|)|a_k|}\rho(\alpha^{n+1}(a_{k})_{\lambda_k}\nonumber\\
  &&\g_{\lambda_i+\lambda_j, \lambda_1,...,\hat{\lambda}_{k, i,j},...,\lambda_{n+2}}([a_{i\lambda_i}a_{j}], \a(a_{1}),... ,\hat{a}_{k, i,j},...,\a(a_{n+2}))\\
  &&+\sum_{1\leq i< j}^{n+2}\sum_{k=i+1}^{j-1}(-1)^{i+j+k+1+A_i+A_j+A_k++|a_i||a_j|+(|\g|+|a_j|)|a_k|}\rho(\alpha^{n+1}(a_{k}))_{\lambda_k}\nonumber\\
  &&\g_{\lambda_i+\lambda_j, \lambda_1,...,\hat{\lambda}_{i, k,j},...,\lambda_{n+2}}([a_{i\lambda_j}a_{j}], \a(a_{1}),... ,\hat{a}_{i, k,j},...,\a(a_{n+2}))\\
  &&+\sum_{1\leq i< j}^{n+2}\sum_{k=j+1}^{n+2}(-1)^{i+j+k+A_i+A_j+A_k+|a_i||a_j|+|\g||a_k|}\rho(\alpha^{n+1}(a_{k})_{\lambda_k}\nonumber\\
  &&\g_{\lambda_i+\lambda_j, \lambda_1,...,\hat{\lambda}_{i, j,k},...,\lambda_{n+2}}([a_{i\lambda_i}a_{j}], \a(a_{1}),... ,\hat{a}_{i,j,k},...,\a(a_{n+2}))\\
  &&+\sum_{1\leq i<j}^{n+2}(-1)^{i+j+A_i+A_j+|a_i||a_j|+|\g|(|a_i|+|a_j|)}\nonumber\\
  &&\rho(\alpha^{n}([a_{i\lambda_i}a_{j})_{\lambda_i+\lambda_j}\g_{\lambda_1,...,\hat{\lambda}_{j},...\hat{\lambda}_{i},
  \lambda_{n+2}}(\a(a_{1}),... ,\hat{a}_{j},...,\hat{a}_{i},...,\a(a_{n+2})\\
  &&+\sum^{n+2}_{distinct i,j,k,l,i< j,k<l}(-1)^{i+j+k+l}sign\{i,j, k, l\}(-1)^{A_i+A_j+|a_i||a_j|+(|a_i|+|a_j|)(|a_k|+|a_l|)}\nonumber\\
  &&\g_{\lambda_k+\lambda_l, \lambda_i+\lambda_j,\lambda_1,...,\hat{\lambda}_{i,j,k,l},...,\lambda_{n+2}}(\a[a_{k\lambda_k}a_{l}],
  \a[a_{i\lambda_i}a_{j}],...,\hat{a}_{i,j,k,l},..., \a^2(a_{n+2}))~~~~~~~~~\\
   &&\sum^{n+2}_{i,j,k=1,i< j,k\neq i,j}(-1)^{i+j+k+l}sign\{i,j, k\}(-1)^{A_i+A_j+|a_i||a_j|+(A_k-|a_i|-|a_j|)}\nonumber\\
  &&\g_{\lambda_i+\lambda_k+\lambda_j,\lambda_1,...,\hat{\lambda}_{i,j,k},...,\lambda_{n+2}}([[a_{i\lambda_i}a_{j}]_{\lambda_i+\lambda_j},\a(\a_k)],
  \a^2(a_1),...,\hat{a}_{i,j,k},..., \a^2(a_{n+2}))~~~~~~~~~
 \end{eqnarray}
where $A_i=(|a_i|+...+|a_{i-1}|)|a_i|$, $A_j-|a_i|=(a_1+...+\hat{a}_i+...+|a_{j-1}|)|a_j|$ and $sign\{i_1, ... , i_p\}$ is the sign of the permutation putting the indices in increasing order and $\hat{a}_{i,j}$,¡¤¡¤¡¤ means that $a_i, a_j,... $ are omitted.

It is obvious that (3.3) and (3.8) summations cancel each other. The same is true
for (3.4) and (3.7), (3.5) and (3.6). The Hom-Jacobi identity implies (3.11) = 0, whereas
skew-symmetry of gives (3.10)=0. As $(M, \b)$ is an $R$-module,
\begin{eqnarray*}
-\rho(\a(a_i))_{\lambda_i}(\rho(a_j)_{\lambda_j} m)+(-1)^{|a_i||a_j|}\rho(\a(a_j))_{\lambda_j}(\rho(a_i)_{\lambda_i} m)=\rho([a_{i\lambda_i}a_j])_{\lambda_i+\lambda_j}\b(m),
\end{eqnarray*}
By $\g\circ \a=\b\circ \g$, we have (3.1),(3.2) and (3.9) summations cancel. This proves $d^2\g=0$. \hfill $\square$

 Thus the cochains of a Hom-Lie conformal superalgebra $R$ with coefficients in a module $M$ form a
complex, which is denoted by
 \begin{eqnarray*}
\widetilde{C}^{\bullet}_{\a} =\widetilde{C}^{\bullet}_{\a}(R,M)=\bigoplus_{n\in \mathbb{Z}_{+}}C^{n}_{\a}(R,M).
 \end{eqnarray*}
This complex is called the basic complex for the $R$-module $(M,\b)$.

 Let $R$ be a regular Hom-Lie conformal superalgebra. Define
 \begin{eqnarray}
\rho(a)_\lambda b=[\a^{s}(a)_\lambda b],~~~\mbox{for any $a,b\in R$}.
 \end{eqnarray}

 Let $\g\in \widetilde{C}^n_{\a}(R,R_s)$.  Define an operator $d_s: \widetilde{C}^{n}_{\a}(R,R_s)\rightarrow \widetilde{C}^{n+1}_{\a}(R,R_s)$ by
\begin{eqnarray*}
&&(d_s\gamma)_{\lambda_1,...,\lambda_{n+1}}(a_{1},..., a_{n+1})\nonumber\\
&=&\sum_{i=1}^{ n+1}(-1)^{i+1}(-1)^{(|\g|+|a_1|+...+|a_{i-1}|)|a_i|}[\a^{n+s}(a_{i})_{\lambda_i}\gamma_{\lambda_1,...,\hat{\lambda}_{i},...,\lambda_{n+1}}(a_{1},..., \hat{a}_{i},..., a_{n+1})]\nonumber\\
  &&+\sum_{1\leq i<j\leq n+1}(-1)^{i+j}(-1)^{(|a_1|+...+|a_{i-1}|)|a_i|+(|a_1|+...+|a_{j-1}|)|a_j|+|a_i||a_j|}\\
   &&\gamma_{\lambda_i+\lambda_j, \lambda_1,...,\hat{\lambda}_i,...\hat{\lambda}_j,...,\lambda_{n+1}}([a_{i\lambda_i}a_{j}], \a(a_{1}),... ,\hat{a}_{i},... , \hat{a}_{j},...,\a(a_{n+1}))~~~~~~~~~
\end{eqnarray*}
 Obviously, the operator $d_s$ is induced from the differential $d$. Thus $d_s$ preserves the space of cochains
and satisfies $d_s^2=0$. In the following the complex $\widetilde{C}^{\bullet}_{\a}(R, R)$ is assumed to be associated with the
differential $d_s$.

Taking $s=-1$, for $\psi\in C^{2}_{\a}(R,R_{-1})_{\overline{0}}$ be a bilinear operator commuting with $\a$,  we consider a $t$-parameterized family of bilinear operations on $R$,
\begin{eqnarray}
[a_\lambda b]_t=[a_\lambda b]+t\psi_{\lambda, -\partial-\lambda}(a,b),~~~\forall a,b\in R.
\end{eqnarray}
If $[\c_{\lambda}\c]$ endows $(R,[\c_{\lambda}\c], \a)$ with a Hom-Lie conformal superalgebra structure, we say that $\psi$  generates
a deformation of the Hom-Lie conformal superalgebra $R$. It is easy to see that $[\c_{\lambda}\c]$ satisfies (2.1) and (2.2).

 If it is true for (2.3), expanding the Hom-Jacobi identity for $[\c_{\lambda}\c]$ gives
 \begin{eqnarray*}
 &&[\a(a)_\lambda[b_{\mu}c]]+t([\a(a)_\lambda(\psi_{\mu,-\partial-\mu}(b,c)]+\psi_{\lambda,-\partial-\lambda}(\a(a),[b_{\mu}c]))\\
 &&+t^2\psi_{\lambda,-\partial-\lambda}(\a(a),\psi_{\mu,-\partial-\mu}(b,c))\\
 &=&(-1)^{|a||b|} [\a(b)_\mu[a_{\lambda}c]]+t([\a(b)_\mu(\psi_{\lambda,-\partial-\lambda}(a,c))]+\psi_{\mu,-\partial-\mu}(\a(b),[a_{\lambda}c]))\\
 &&+(-1)^{|a||b|}t^2\psi_{\mu,-\partial-\mu}(\a(b),\psi_{\lambda,-\partial-\lambda}(a,c))+[[a_{\lambda}b]_{\lambda+\mu}\a(c)]\\
 &&+t([(\psi_{\lambda,-\partial-\lambda}(a,b))_{\lambda+\mu}\a(c)]+ \psi_{\lambda+\mu,-\partial-\lambda-\mu}([a_\lambda b],\a(c)))\\
 &&+t^2\psi_{\lambda+\mu,-\partial-\lambda-\mu}(\psi_{\lambda,-\partial-\lambda}(a,b), \a(c)).
 \end{eqnarray*}
This is equivalent to the following conditions
\begin{eqnarray}
 &&\psi_{\lambda,-\partial-\lambda}(\a(a),\psi_{\mu,-\partial-\mu}(b,c))\nonumber\\
 &=&(-1)^{|a||b|}\psi_{\mu,-\partial-\mu}(\a(b),\psi_{\lambda,-\partial-\lambda}(a,c))+ \psi_{\lambda+\mu,-\partial-\lambda-\mu}(\psi_{\lambda,-\partial-\lambda}(a,b),\a(c)).~~~~~
\end{eqnarray}
and
\begin{eqnarray}
&&[\a(a)_\lambda(\psi_{\mu,-\partial-\mu}(b,c)]+\psi_{\lambda,-\partial-\lambda}(\a(a),[b_{\mu}c])\nonumber\\
&=& (-1)^{|a||b|}[\a(b)_\mu(\psi_{\lambda,-\partial-\lambda}(a,c))]+(-1)^{|a||b|}\psi_{\mu,-\partial-\mu}(\a(b),[a_{\lambda}c])\nonumber\\
&&+[(\psi_{\lambda,-\partial-\lambda}(a,b))_{\lambda+\mu}\a(c)]+ \psi_{\lambda+\mu,-\partial-\lambda-\mu}([a_\lambda b],\a(c)).
\end{eqnarray}

By conformal antilinearity of $\psi$, we have
\begin{eqnarray}
[(\psi_{\lambda, -\partial-\lambda}(a,b))_{\lambda+\mu}c]=[\psi_{\lambda, \mu}(a,b)_{\lambda+\mu}c].
\end{eqnarray}

On the other hand, Let $\psi$ be a cocycle, i.e., $d_{-1}\psi=0$. In fact,
\begin{eqnarray}
0&=&(d^{-1}\psi)_{\lambda, \mu, \g}(a,  b, c)\nonumber\\
&=& [\a(a)_{\lambda}\psi_{\mu, \g}(b, c)]-(-1)^{|a||b|}[\a(b)_{\mu}\psi_{\lambda, \g}(a, c)]+(-1)^{(|a|+|b|)|c|}[\a(c)_{\g}\psi_{\lambda, \mu}(a, b)]\nonumber\\
&& -\psi_{\lambda+\mu, \g}([a_\lambda b], \a(c))+(-1)^{|b||c|}\psi_{\lambda+\g, \mu}([a_\lambda c], \a(b))-(-1)^{|a|(|b|+|c|)}\psi_{\g+\mu, \lambda}([b_\mu c], \a(a))\nonumber\\
&=& [\a(a)_{\lambda}\psi_{\mu, \g}(b, c)]-(-1)^{|a||b|}[\a(b)_{\mu}\psi_{\lambda, \g}(a, c)]-\psi_{\lambda, \mu}(a, b)_{-\partial-\g}\a(c)]\nonumber\\
&& +\psi_{\lambda, \g+\mu}(\a(a), [b_\mu c])-(-1)^{|b||a|}\psi_{\mu, \lambda+\g}(\a(b), [a_\lambda c])-\psi_{\lambda+\mu, \g}([a_\lambda b], \a(c)).
\end{eqnarray}

Replacing $\g$ by $-\lambda-\mu-\partial$ in (3.17), we have
\begin{eqnarray*}
&&[\a(a)_\lambda(\psi_{\mu,-\partial-\mu}(b,c)]+\psi_{\lambda,-\partial-\lambda}(\a(a),[b_{\mu}c])\nonumber\\
&=& (-1)^{|a||b|}[\a(b)_\mu(\psi_{\lambda,-\partial-\lambda}(a,c))]+(-1)^{|a||b|}\psi_{\mu,-\partial-\mu}(\a(b),[a_{\lambda}c])\nonumber\\
&&+[(\psi_{\lambda,-\partial-\lambda}(a,b))_{\lambda+\mu}\a(c)]+ \psi_{\lambda+\mu,-\partial-\lambda-\mu}([a_\lambda b],\a(c)).
\end{eqnarray*}

When $\psi$ is a 2-cocycle satisfying (3.14), $(R,[\c_{\lambda}\c], \a)$ forms  a Hom-Lie conformal superalgebra. In this case, $\psi$ generates a deformation of the Lie conformal superalgebra $R$.

A deformation is said to be trivial if there is a linear operator $f\in \widetilde{C}_{\a}^{1}(R,R_{-1})$ such that for $T_t=id+tf$, there holds
\begin{eqnarray}
T_t([a_\lambda b]_t)=[T_t(a)_\lambda T_t(b)],~~~\forall a.b\in R.
\end{eqnarray}
\begin{definition}
 A linear operator $f\in \widetilde{C}_{\a,\b}^{1}(R,R_{-1})$ is a Hom-Nijienhuis operator
if  \begin{eqnarray}
    [f(a)_\lambda f(b)]=f([a_\lambda b]_{N}),~~~\forall a.b\in R.
    \end{eqnarray}
where the bracket $[\c,\c]_{N}$ is defined by
\begin{eqnarray}
[a_\lambda b]_{N}=[f(a)_\lambda b]+[a_\lambda f(b)]-f([a_\lambda b]),~~~\forall a.b\in R.
\end{eqnarray}
\end{definition}
\begin{theorem}
 Let $(R,\a)$ be a regular Hom-Lie conformal superalgebra, and $f\in \widetilde{C}_{\a}^{1}(R,R_{-1})$
a Hom-Nijienhuis operator. Then a deformation of $(R,\a)$) can be obtained by putting
\begin{eqnarray}
\psi_{\lambda, -\partial-\lambda}(a,b)=[a_\lambda b]_{N},~~~\forall a, b\in R.
\end{eqnarray}
Furthermore, this deformation is trivial.
\end{theorem}
{\bf Proof.}  To see that $\psi$  generates a deformation, we need to check   satisfying (3.14), by (3.20)
and (3.21), we have
\begin{eqnarray*}
&&\psi_{\lambda,-\partial-\lambda}(\a(a),\psi_{\mu,-\partial-\mu}(b,c))-(-1)^{|a||b|}\psi_{\mu,-\partial-\mu}(\a(b),\psi_{-\partial-\lambda,\lambda}(c,a))\\
&&-\psi_{-\partial-\lambda-\mu, \lambda+\mu}(\a(c), \psi_{\lambda,-\partial-\lambda}(a,b))\\
&=&[f(\a(a))_{\lambda}[f(b)_{\mu}c]]+[f(\a(a))_{\lambda}[b_\mu f(c)]]-[f(\a(a))_{\lambda}f([b_{\mu}c])] \\
&&+[\a(a)_{\lambda}[f(b)_\mu f(c)]]-f([\a(a)_{\lambda}[f(b)_\mu c]])-f([\a(a)_{\lambda}[b_\mu f(c)]])+f([\a(a)_{\lambda}f([b_\mu c])])\\
&&-(-1)^{|a||b|}[f(\a(b))_{\mu}[f(c)_{-\partial-\lambda}a]]-(-1)^{|a||b|}[f(\a(b))_{\mu}[c_{-\partial-\lambda} f(a)]]\\
&&+(-1)^{|a||b|}[f(\a(b))_{\mu}f([c_{-\partial-\lambda}a])]-(-1)^{|a||b|}[\a(b)_{\mu}[f(c)_{-\partial-\lambda} f(a)]]\\
&&+(-1)^{|a||b|}f([\a(b)_{\mu}[f(c)_{-\partial-\lambda} a]])+(-1)^{|a||b|}f([\a(b)_{\mu}[c_{-\partial-\lambda} f(a)]])\\
&&-(-1)^{|a||b|}[f(\a(b))_{\mu}f([c_{-\partial-\lambda} a])])-[f(\a(c))_{-\partial-\lambda-\mu}[f(a)_{\lambda}b]]-[f(\a(c))_{-\partial-\lambda-\mu}[a_\lambda f(b)]]\\
&&+[f(\a(c))_{-\partial-\lambda-\mu}f([a_{\lambda}b])]-[\a(c)_{-\partial-\lambda-\mu}[f(a)_{\lambda}f(b)]]\\
&&+f([\a(c)_{-\partial-\lambda-\mu}[f(a)_{\lambda}b]])+f([\a(c)_{-\partial-\lambda-\mu}[a_{\lambda}f(b)]])-f([\a(c)_{-\partial-\lambda-\mu}f([a_{\lambda}b])]).
\end{eqnarray*}
Since $f$ is a Hom-Nijenhuis operator, we get
\begin{eqnarray*}
&&-[f(\a(a))_{\lambda}f([b_{\mu}c])]+f([\a(a)_{\lambda}f([b_\mu c])])\\
&=&-f([f(\a(a))_{\lambda}[b_{\mu}c]])+f^2([\a(a)_{\lambda}[b_\mu c]]),
\end{eqnarray*}
\begin{eqnarray*}
&& -(-1)^{|a||b|}[f(\a(b))_{\mu}f([a_{\lambda}c])]+(-1)^{|a||b|}f([\a(b)_{\mu}f([a_{\lambda} c])])\\
&=&-(-1)^{|a||b|}f([f(\a(b))_{\mu}[a_{\lambda}c]])+(-1)^{|a||b|}f^2([\a(b)_{\mu}[a_{\lambda} c]])
\end{eqnarray*}
and
\begin{eqnarray*}
&& -[f(\a(c))_{-\partial-\lambda-\mu}f([a_{\lambda}b])]+f([c_{-\partial-\lambda-\mu}f([a_{\lambda}b])])\\
&=&-f[f(\a(c))_{-\partial-\lambda-\mu}[a_{\lambda}b]]+f^2([c_{-\partial-\lambda-\mu}[a_{\lambda}b]])
\end{eqnarray*}
Note that
\begin{eqnarray*}
[\a(a)_\lambda[f(b)_\mu f(c)]]=[[a_{\lambda}f(b)]_{\lambda+\mu}f(\a(c))]+(-1)^{|a||b|}[f(\a(b))_{\mu}[a_{\lambda} f(c)]].
\end{eqnarray*}

Thus
\begin{eqnarray*}
 &&\psi_{\lambda,-\partial-\lambda}(\a(a),\psi_{\mu,-\partial-\mu}(b,c))\nonumber\\
 &=&(-1)^{|a||b|}\psi_{\mu,-\partial-\mu}(\a(b),\psi_{\lambda,-\partial-\lambda}(a,c))+ \psi_{\lambda+\mu,-\partial-\lambda-\mu}(\psi_{\lambda,-\partial-\lambda}(a,b),\a(c)).
\end{eqnarray*}
In the similar way, we can check that  (3.15). This proves that $\psi$ generates a deformation of the regular Hom-Lie conformal superalgebra $(R, \a)$.

Let $T_t=id+tf$. By  (3.13) and (3.21), we have
\begin{eqnarray}
T_t([a_\lambda b]_t)&=&(id+tf)([a_\lambda b]+t\psi_{\lambda,-\partial-\lambda}(a,b))\nonumber\\
&=& (id+tf)([a_\lambda b]+t[a_\lambda b]_N)\nonumber\\
&=& [a_\lambda b]+t([a_\lambda b]_N+f([a_\lambda b]))+t^2f([a_\lambda b]_N).
\end{eqnarray}

On the other hand, we have
\begin{eqnarray}
[T_t(a)_\lambda T_t(b)]&=&[(a+tf(a))_\lambda(b+tf(b))]\nonumber\\
&=& [a_\lambda b]+t([f(a)_\lambda b]+[a_\lambda f(b)])+t^2[f(a)_\lambda f(b)].
\end{eqnarray}
Combining  (3.22) with (3.23) gives $T_t([a_\lambda b]_t)=[T_t(a)_\lambda T_t(b)]$. Therefore the deformation
is trivial.\hfill $\square$

\section{Derivations of multiplicative  Hom-Lie conformal superalgebras }
\def\theequation{\arabic{section}. \arabic{equation}}
\setcounter{equation} {0}

\begin{definition} Let $(R,\a)$ be a multiplicative Hom-Lie conformal superalgebra. Then a
Hom-conformal linear map $D_\lambda: R\rightarrow R$ is called an $\a^k$-derivation of $(R,\a)$ if
\begin{eqnarray}
D_\lambda \circ \a=\a\circ D_\lambda,\nonumber\\
D_\lambda([a_\mu b])=[D_\lambda(a)_{\lambda+\mu}\a^{k}(b)]+(-1)^{|a||D|}[\a^{k}(a)_{\mu}D_\lambda(b)].
\end{eqnarray}
\end{definition}
Denote by $Der_{\a^{s}}$ the set of $\a^{s}$-derivations of the multiplicative Hom-Lie conformal
superalgebra $(R,\a)$. For any $a\in R$ satisfying $\a(a)=a$, define $D_{k}:R\rightarrow R$ by
\begin{eqnarray*}
D_{k}(a)_{\lambda}(b)=[a_\lambda\a^{k+1}(b)],~~~\forall b\in R.
\end{eqnarray*}

Then $D_{k}(a)$ is an $\a^{k+1}$-derivation, which is called an inner $\a^{k+1}$-derivation. In fact
\begin{eqnarray*}
D_{k}(a)_{\lambda}(\partial b)&=&[a_\lambda\a^{k+1}(\partial b)]\\
&=& [a_\lambda\partial\a^{k+1}( b)]\\
&=& (\partial+\lambda)D_{k}(a)_{\lambda}(b),\\
D_{k}(a)_{\lambda}(\a(b))&=&[a_\lambda\a^{k+2}(b)]\\
&=& \a[a_\lambda\a^{k+1}( b)]\\
&=& \a\circ D_{k}(a)_{\lambda}(b),\\
D_{k}(a)_{\lambda}([b_{\mu}c])&=& [a_\lambda \a^{k+1}([b_{\mu}c])\\
&=&  [\a(a)_\lambda [\a^{k+1}(b)_{\mu}\a^{k+1}(c)]\\
&=& [a_\lambda\a^{k+1}(b)]_{\lambda+\mu}\a^{k+1}(c)]+(-1)^{|a||b|}[\a^{k+1}(b)_{\mu}[\a(a)_{\lambda}\a^{k+1}(c)]]\\
&=&[a_\lambda\a^{k+1}(b)]_{\lambda+\mu}\a^{k+1}(c)]+(-1)^{|a||b|}[\a^{k+1}(b)_{\mu}[a_{\lambda}\a^{k+1}(c)]]\\
&=&[D_{k}(a)_{\lambda}(b)_{\lambda+\mu}\a^{k+1}(c)]+(-1)^{|a||b|}[\a^{k+1}(b)_{\mu}(D_{k}(a)_{\lambda}(c))]
\end{eqnarray*}
Denote by $Inn_{\a^k} (R)$ the set of inner $\a^k$-derivations. For $D_{\lambda}\in Der_{\a^k}(R)$ and $D'_{\mu-\lambda}\in Der_{\a^s}(R)$, define their commutator $[D_\lambda D']_{\mu}$ by
\begin{eqnarray}
[D_\lambda D']_{\mu}(a)=D_\lambda(D'_{\mu-\lambda}a)-(-1)^{|D||D'|}D'_{\mu-\lambda}(D_\lambda a),~~~\forall a\in R.
\end{eqnarray}

\begin{lemma} For any $D_{\lambda}\in Der_{\a^k}(R)$ and $D'_{\mu-\lambda}\in Der_{\a^s}(R)$, we have
\begin{eqnarray*}
[D_\lambda D'] \in Der_{\a^{k+s}}(R)[\lambda].
\end{eqnarray*}
\end{lemma}

{\bf Proof.} For any $a,b\in R$, we have
\begin{eqnarray*}
&&[D_\lambda D']_{\mu}(\partial a)\\
&=&D_{\lambda}(D'_{\mu-\lambda}\partial a)-(-1)^{|D||D'|}D'_{\mu-\lambda}(D_\lambda \partial a)\\
&=& D_\lambda((\partial+\mu-\lambda)D'_{\mu-\lambda}a)+(-1)^{|D||D'|}D'_{\mu-\lambda}((\mu+\lambda)D_\lambda a)\\
&=& (\partial+\mu)D_{\lambda}(D'_{\mu-\lambda}a)-(-1)^{|D||D'|}(\partial+\mu)D'_{\mu-\lambda}(D_{\lambda}a)\\
&=& (\partial+\mu)[D_\lambda D']_{\mu}(a).
\end{eqnarray*}
and
\begin{eqnarray*}
&&[D_\lambda D']_{\mu}([a_{\g}b])\\
&=& D_\lambda(D'_{\mu-\lambda}[a_{\g}b])-(-1)^{|D||D'|}D'_{\mu-\lambda}(D_\lambda[a_{\g}b])\\
&=&D_\lambda([D'_{\mu-\lambda}(a)_{\mu-\lambda+\gamma}\a^s(b)]+(-1)^{|x||D|}[\a^s(a)_{\g}D'_{\mu-\lambda}(b)])\\
&& -(-1)^{|D||D'|}D'_{\mu-\lambda}([D_\lambda(a)_{\lambda+\gamma}\a^k(b)]+(-1)^{|x||D|}[\a^k(a)_{\g}D_{\lambda}(b)])\\
&=& [D_\lambda(D'_{\mu-\lambda}(a))_{\mu+\gamma}\a^{k+s}(b)]+(-1)^{|D||D'(x)|}[\a^{k}(D'_{\mu-\lambda}(a))_{\mu-\lambda+\g}D_{\lambda}(\a^{s}(b)))\\
&&+(-1)^{|x||D'|}[D_\lambda(\a^s(a))_{\lambda+\gamma}\a^k(D'_{\mu-\lambda}(b))]+(-1)^{|x||D|}[\a^{k+s}(a)_\g(D_\lambda(D'_{\mu-\lambda}(b)))]\\
&&-(-1)^{|D||D'|}[(D'_{\mu-\lambda}D_\lambda(a))_{\mu+\gamma}\a^{k+s}(b)]+(-1)^{|D'||D(x)|}[\a^{s}(D_\lambda(s))_{\lambda+\gamma}(D'_{\mu-\lambda}(\a^k(b)))]\\
&&-(-1)^{|D|(|D'|+|x|)}[(D'_{\mu-\lambda}(\a^{k}(a)))_{\mu-\lambda+\gamma}\a^{s}(D_{\lambda}(b))]+(-1)^{|x||D'|}[\a^{k+s}(a)_{\lambda}(D'_{\mu-\lambda}(D_\lambda(b)))]\\
&=&[([D_\lambda D']_\mu a)_{\mu+\g}\a^{k+s}(b)]+(-1)^{|x||[D,D']|}[\a^{k+s}(a)_{\g}([D_\lambda D']_{\mu}b)].
\end{eqnarray*}
Therefore, $[D_\lambda D'] \in Der_{\a^{k+s}}(R)[\lambda]$. \hfill $\square$

Define
\begin{eqnarray}
Der(R)=\bigoplus_{k\geq0, l\geq 0}Der_{\a^k}(R).
\end{eqnarray}

\begin{proposition}
$(Der(R), \a')$ is a Hom-Lie conformal superalgebra with respect to (4.2), where $\a'(D)=D\circ \a$.
\end{proposition}

{\bf Proof.}  By  (4.2), it is easy to check that (2.1) and (2.2)
are satisfied. To check the Hom-Jacobi identity, we compute separately
\begin{eqnarray*}
&&[\a'(D)_{\lambda}[D'_{\mu}D'']]_{\theta}(a)\\
&=&(D\circ \a)_{\lambda}([D'_{\mu}D'']_{\theta-\lambda}a)-[D'_{\mu}D'']_{\theta-\lambda}((D\circ \a)_\lambda a)\\
&=& D_{\lambda}([D'_\mu D'']_{\theta-\lambda}\a(a))-(-1)^{|D|(|D'|+|D''|)}[D'_\mu D'']_{\theta-\lambda}(D_\lambda\a(a))\\
&=& D_\lambda(D'_\mu(D''_{\theta-\lambda-\mu}\a(a)))-(-1)^{|D''||D'|)}D_\lambda(D''_{\theta-\lambda-\mu}(D'_{\mu}\a(a)))\\
&& -(-1)^{|D||D'|}D_{\mu}'(D''_{\theta-\lambda-\mu}(D_\lambda\a(a)))+(-1)^{|D||D'|+|D||D''|+|D'||D''|}D''_{\theta-\lambda-\mu}(D'_\mu(D_\lambda\a(a))),\\
&&(-1)^{|D||D'|}[\a'(D')_{\mu}[D_\lambda D'']]_{\theta}(a)\\
&=& (-1)^{|D||D'|}D'_{\mu}(D_\lambda(D''_{\theta-\lambda-\mu}\a(a)))-(-1)^{|D||D'|}D'_\mu(D''_{\theta-\lambda-\mu}(D_\lambda(\a(a))))\\
&& -(-1)^{|D'||D''|}D_\lambda(D''_{\theta-\lambda-\mu}(D'_\mu\a(a)))+(-1)^{(|D|+|D'|)|D''|}D''_{\theta-\lambda-\mu}(D_\lambda(D'_\mu\a(a))),\\
&&[[D_\lambda D']_{\lambda+\mu}\a'(D'')]_{\theta}(a)\\
&=&D_\lambda(D'_\mu(D''_{\theta-\lambda-\mu}\a(a)))-(-1)^{|D||D'|}D'_{\mu}(D_\lambda(D''_{\theta-\lambda-\mu}\a(a)))\\
&&-(-1)^{(|D|+|D'|)|D''|}D''_{\theta-\lambda-\mu}(D_\lambda(D'_\mu\a(a)))+(-1)^{|D||D'|+|D||D''|+|D'||D''|}D''_{\theta-\lambda-\mu}(D'_\mu(D_\lambda\a(a))).
\end{eqnarray*}
Thus, $[\a'(D)_{\lambda}[D'_{\mu}D'']]_{\theta}(a)=(-1)^{|D||D'|}[\a'(D')_{\mu}[D_\lambda D'']]_{\theta}(a)+[[D_\lambda D']_{\lambda+\mu}\a'(D'')]_{\theta}(a)$.  This proves that $(Der(R), \a')$ is a Hom-Lie conformal superalgebra. \hfill $\square$

At the end of this section, we give an application of the $\a$-derivations of a  regular
Hom-Lie conformal superalgebra $(R, \a)$. For any $D_\lambda\in Cend(R)$,  define a bilinear
operation $[\c_\lambda \c]_D$ on the vector space $R\oplus \mathbb{R}D$ by
\begin{eqnarray}
[(a+mD)_\lambda(b+nD)]_D=[a_\lambda b]+mD_\lambda(b)-(-1)^{|a||b|}nD_{-\lambda-\partial}(a), \forall a,b\in R, m,n\in \mathbb{R}.~~~~~~
\end{eqnarray}
and a linear map $\a': R\oplus \mathbb{R}D\rightarrow R\oplus \mathbb{R}D$ by $\a'(a+D)=\a(a)+D$.

\begin{proposition}
  $(R\oplus \mathbb{R}D, \a')$ is a  regular Hom-Lie conformal superalgebra if and only
if $D_\lambda$ is an $\a$-derivation of $(R, \a)$.
\end{proposition}

{\bf Proof.} Suppose that $(R\oplus \mathbb{R}D, \a')$ is a   regular Hom-Lie conformal superalgebra. For any $a,b\in R$, $m,n\in \mathbb{R}$, we have
\begin{eqnarray*}
&&\a'[(a+mD)_\lambda(b+nD)]_D\\
&=&\a'([a_\lambda b]+mD_\lambda(b)-(-1)^{|a||b|}nD_{-\lambda-\partial}(a))\\
&=& \a[a_\lambda b]+m\a (D_\lambda(b))-(-1)^{|a||b|}n\a(D_{-\lambda-\partial}(a)),
\end{eqnarray*}
and
\begin{eqnarray*}
&&[\a'(a+mD)_{\lambda}\a'(b+nD)]\\
&=&[\a(a)+mD_{\lambda}\a(b)+nD]\\
&=&[\a(a)_\lambda \a(b)]+mD_\lambda \a(b)-(-1)^{|a||b|}nD_{-\lambda-\partial} \a(a).
\end{eqnarray*}
Thus, we have
\begin{eqnarray*}
\a\circ D_\lambda=D_\lambda\circ \a.
\end{eqnarray*}
Next, the Hom Jacobi identity gives
\begin{eqnarray*}
 [\a'(D)_{\lambda}[a_{\mu}b]_D]_D=(-1)^{|a||b|}[\a'(a)_{\lambda}[D_\mu b]_D]_{D}+[[D_\mu a]_{D\lambda+\mu}\a'(b)]_{D}
\end{eqnarray*}
which is exactly $D_{\mu}([a_\lambda b])=[(D_\mu a)_{\lambda+\mu}\a(b)]+(-1)^{|a||b|}[\a(a)_{\lambda}(D_\mu b)]$ by (4.4). Therefore, $D_\lambda$ is an $\a$-derivation of $(R, \a)$.

Conversely, let $D_\lambda$ is an $\a$-derivation of $(R, \a)$. For any $a,b\in R$, $m,n\in \mathbb{R}$,
\begin{eqnarray*}
&&[(b+nD)_{-\partial-\lambda}(a+mD)]_D\\
&=& [b_{-\partial-\lambda}a]+nD_{-\partial-\lambda}(a)-(-1)^{|a||b|}mD_\lambda(b)\\
&=&-(-1)^{|a||b|}([a_\lambda b]-(-1)^{|a||b|}nD_{-\partial-\lambda}(a)+mD_\lambda(b))\\
&=& -(-1)^{|a||b|}[(a+mD)_{\lambda}(b+nD)]_D
\end{eqnarray*}
which proves (2.2). And it is obvious that
\begin{eqnarray*}
&&[\partial D_\lambda a]_D=-\lambda[D_\lambda a]_D,\\
&& [\partial a_\lambda D]_D=-D_{-\partial-\lambda}(\partial a)=-\lambda[a_\lambda D]_D,\\
&& [D_\lambda \partial a]_D=D_\lambda(\partial a)=(\partial+\lambda)D_\lambda(a)=(\partial+\lambda)[D_\lambda a]_D,\\
&&[a_\lambda \partial D]_D=-(\partial D)_{-\lambda-\partial}a=(\partial+\lambda)[a_\lambda D]_D,\\
&& \a'\circ \partial =\partial\circ \a'.
\end{eqnarray*}
Thus (2.1) follows. The Hom Jacobi identity is easy to check.\hfill $\square$

\section{Generalized derivations of multiplicative  Hom-Lie conformal superalgebras }
\def\theequation{\arabic{section}. \arabic{equation}}
\setcounter{equation} {0}

Let $(R,\a)$ be a  multiplicative Hom-Lie conformal superalgebra. Define
\begin{eqnarray*}
\Omega=\{D_\lambda\in Cend(R)|D_\lambda \circ \a=\a\circ D_\lambda\}.
\end{eqnarray*}
Then $\Omega$ is a Hom-Lie conformal superalgebra, and $Der(R)$ is a subalgebra of $\Omega$.

\begin{definition}  An element $D_{\mu}$ in $\Omega$ is called\\
(a) a generalized derivation  of $R$, if there exist $D'_{\mu}, D''_{\mu}\in \Omega$, such that
\begin{eqnarray}
[(D_\mu(a))_{\lambda+\mu}\a^{k}(b)]+(-1)^{|D||a|}[\a^{k}(a)_{\lambda}(D'_{\mu}(b))]=D''_{\mu}([a_\lambda b]),~~\forall a,b\in R.
\end{eqnarray}
(b) an $\a^k$-quasiderivation  of $R$, if there is $D'_{\mu}\in \Omega$ such that
\begin{eqnarray}
[(D_\mu(a))_{\lambda+\mu}\a^{k}(b)]+(-1)^{|D||a|}[\a^{k}(a)_{\lambda}(D_{\mu}(b))]=D'_{\mu}([a_\lambda b]),~~\forall a,b\in R.
\end{eqnarray}
(c) an $\a^k$-centroid  of $R$, if it satisfies
\begin{eqnarray}
[(D_\mu(a))_{\lambda+\mu}\a^{k}(b)]=(-1)^{|D||a|}[\a^{k}(a)_{\lambda}(D_{\mu}(b))]=D_{\mu}([a_\lambda b]),~~\forall a,b\in R.
\end{eqnarray}
(d) an $\a^k$-quasicentroid of $R$, if it satisfies
\begin{eqnarray}
[(D_\mu(a))_{\lambda+\mu}\a^{k}(b)]=(-1)^{|D||a|}[\a^{k}(a)_{\lambda}(D_{\mu}(b))], ~~\forall a,b\in R.
\end{eqnarray}
(e) an $\a^k$-central derivation  of $R$, if it satisfies
\begin{eqnarray}
[(D_\mu(a))_{\lambda+\mu}\a^{k}(b)]=D_{\mu}([a_\lambda b])=0, ~~\forall a,b\in R.
\end{eqnarray}
\end{definition}

Denote by $GDer_{\a^k}(R), QDer_{\a^k}(R), C_{\a^k}(R), QC_{\a^k}(R)$ and $ZDer_{\a^k}(R)$ the sets of
all generalized $\a^k$-derivations, $\a^k$-quasiderivations, $\a^k$-centroids, $\a^k$-quasicentroids and
$\a^k$-central derivations  of $R$, respectively. Set
\begin{eqnarray*}
&&GDer(R):=\bigoplus_{k\geq0}GDer_{\a^k}(R),~~QDer(R):=\bigoplus_{k\geq0}QDer_{\a^k}(R).\\
&& C(R):= \bigoplus_{k\geq0} C_{\a^k}(R),~~~QC_{\a^k}(R):=  \bigoplus_{k\geq0} QC_{\a^k}(R),\\
&&ZDer(R):=\bigoplus_{k\geq0}ZDer_{\a^k}(R).
\end{eqnarray*}
It is easy to see that
\begin{eqnarray}
ZDer(R)\subseteq Der(R)\subseteq QDer(R)\subseteq GDer(R)\subseteq Cend(R),~~ C(R)\subseteq QC(R)\subseteq GDer(R).~~~~~
\end{eqnarray}

\begin{proposition} Let $(R,\a)$ be a  multiplicative Hom-Lie conformal superalgebra. Then

(1) $GDer(R), QDer(R)$  and $C(R)$ are subalgebras of $\Omega$,

(2) $ZDer(R)$ is an ideal of $Der(R)$.
\end{proposition}

{\bf Proof.} (1)We only prove that $GDer(R)$ is a subalgebra of $\Omega$.  The proof for the other two cases is left to reader.

For any $D_{\mu}\in GDer_{\a^k}(R),H_{\mu}\in GDer_{\a^s}(R),a,b\in R$, there exist $D'_{\mu},D''_{\mu}\in \Omega,(resp.H'_{\mu},H''_{\mu}\in \Omega)$
such that (5.1) holds for  $D_{\mu}(resp. H_{\mu})$. Recall that $\alpha'(D_{\mu})=D_{\mu}\circ\alpha$.
\begin{eqnarray*}
&&[(\a'(D_{\mu})(a))_{\lambda+\mu}\alpha^{k+1}(b)]\\
&=&[(D_{\mu}(\alpha(a)))_{\lambda+\mu}\alpha^{k+1}(b)]=\a([(D_{\mu}(a))_{\lambda+\mu}\alpha^{k}(b)])\\
&=&\a(D''_{\mu}([a_{\lambda}b])-(-1)^{|D||a|}[\a^{k}(a)_{\lambda}D'_{\mu}(b)])\\
&=&\a'(D''_{\mu})([a_{\lambda}b])-(-1)^{|D||a|}[\a^{k+1}(a)_{\lambda}(\alpha'(D'_{\mu})(b))].
\end{eqnarray*}
This gives $\alpha'(D_{\mu})\in GDer_{\a^{k+1}}(R)$.  Furthermore, we need to show
\begin{eqnarray}
[D''_{\mu}H'']_{\theta}([a_{\lambda}b])=[([D_{\mu}H]_{\theta}(a))_{\lambda+\theta}\alpha^{k+s}(b)]+(-1)^{(|D''|+|H''|)|a|}[\alpha^{k+s}(a)_{\lambda}([D'_{\mu}H']_{\theta}(b))].~~~~
\end{eqnarray}
By (5.4), we have
\begin{eqnarray}
[([D_{\mu}H]_{\theta}(a))_{\lambda+\theta}(b)]=[(D_{\mu}(H_{\theta-\mu}(a)))_{\lambda+\theta}\alpha^{k+s}(b)]-(-1)^{|D||H|}
[(H_{\theta-\mu}(D_{\mu}(a))))_{\lambda+\theta}\alpha^{k+s}(b)]~~~.
\end{eqnarray}
By (5.1), we obtain
\begin{eqnarray}
&&[(D_{\mu}(H_{\theta-\mu}(a)))_{\lambda+\theta}\alpha^{k+s}(b)]\nonumber\\
&=&D''_{\mu}([(H_{\theta-\mu}(a))_{\lambda+\theta-\mu}\alpha^{s}(b)])-(-1)^{|D|(|H|+|a|)}D''_{\mu}([\alpha^{k}(a)_{\lambda}(H'_{\theta-\mu}(b)))]\nonumber\\
&=&D''_{\mu}(H''_{\theta-\mu}([a_{\lambda}b]))-(-1)^{|H||a|}D''_{\mu}([\alpha^{s}(a)_{\lambda}(H'_{\theta-\mu}(b)))]\\
&&-(-1)^{|D|(|H|+|a|)}H''_{\theta-\mu}([\alpha^{k}(a)_{\lambda}(D'_{\mu}(b))])+(-1)^{|D|(|H|+|a|)+|H||a|}[\alpha^{k+s}(a)_{\lambda}(H'_{\theta-\mu}(D'_{\mu}(b)))]\nonumber\\
&&[(H_{\theta-\mu}(D_{\mu}(a)))_{\lambda+\theta}\alpha^{k+s}(b)]\nonumber\\
&=&H''_{\theta-\mu}([(D_{\mu}(a))_{\lambda+\mu}\alpha^{k}(b)])-(-1)^{|H|(|D|+|a|)}[\alpha^{s}(D_{\mu}(a))_{\lambda+\mu}(H'_{\theta-\mu}(\alpha^{k}(b)))]\nonumber\\
&=&H''_{\theta-\mu}(D''_{\mu}([a_{\lambda}b]))-(-1)^{|D||a|}(H''_{\theta-\mu}([\alpha^{k}(a)_{\lambda}(D'_{\mu}(b)])\nonumber\\
&&-(-1)^{|H|(|D|+|a|)}D''_{\mu}([\alpha^{s}(a)_{\lambda}(H'_{\theta-\mu}(b))])+(-1)^{|H|(|D|+|a|+|D||a|}\nonumber\\
&&[\alpha^{k+s}(a)_{\lambda}(D'_{\mu}(H'_{\theta-\mu}(b))).~~~~
\end{eqnarray}

  Substituting (5.9) and (5.10) into (5.8) gives (5.7). Hence $[D_{\mu}H)\in GDer_{\alpha^{k+s}}(R)[\mu]$, and  $GDer(R)$ is a Hom sub-algebra of
$\Omega$.

(2) For $D_1{\mu}\in ZDer_{\a^k}(R),D_2{\mu}\in Der_{\a^s}(R)$, and  $a,b\in R$, we have
\begin{eqnarray*}
[(\alpha'(D_1)_{\mu}(a)_{\lambda+\mu}\alpha^{k+1}(b)]=\alpha([(D_1{\mu}(a))_{\lambda+\mu}\alpha^{k}(b)])=\alpha'(D_1)_{\mu}([a_{\lambda}b])=0,
\end{eqnarray*}
which proves $\alpha'(D_1)\in ZDer_{\a^{k+1}}(R).$   By (5.5),
\begin{eqnarray*}
&&[D_1{\mu}D_{2}]_{\theta}([a_{\lambda}b])\\
&=&D_1{\mu}(D_{2\theta-\mu}([a_{\lambda}b]))-(-1)^{|D_1||D_2|}D_{2\theta-\mu}(D_1{\mu}([a_{\lambda}b]))\\
&=&D_1{\mu}(D_{2\theta-\mu}([a_{\lambda}b]))\\
&=&D_1{\mu}([(D_{2\theta-\mu}(a))_{\lambda+\theta-\mu}\alpha^{s}(b)]+(-1)^{|D_2||a|}[\alpha^{s}(a)_{\lambda}D_{2\theta-\mu}(b)])=0,\\
&&[D_1{\mu}D_{2}]_{\theta}(a)_{\lambda+\theta}\alpha^{k+s}(b)]\\
&=&[D_1{\mu}(D_{2\theta-\mu}(a))-(-1)^{|D_1||D_2|}D_{2\theta-\mu}(D_1{\mu}(a)))_{\lambda+\theta}\alpha^{k+s}(b)]\\
&=&[-(-1)^{|D_1||D_2|}(D_{2\theta-\mu}(D_1{\mu}(a)))_{\lambda+\theta}\alpha^{(k+s)}(b)]\\
&=&-(-1)^{|D_1||D_2|}(D_{2\theta-\mu}([D_1{\mu}(a)_{\lambda+\mu}\alpha^{k}(b)])+(-1)^{|D_2||a|}\alpha^{s}(D_1{\mu}(a))_{\lambda+\mu}D_{2\theta-\mu}(\alpha^{k}(b))]\\
&=&0.
\end{eqnarray*}
This shows that $[D_1{\mu}D_{2}]\in ZDer_{\a^{k+s}}(R)[\mu].$  Thus $ZDer(R)$ is an ideal of $Der(R)$. \hfill $\square$

\begin{lemma}
 Let $(R,\a)$ be a  multiplicative Hom-Lie conformal superalgebra. Then \\
(1) $[Der(R)_\lambda C(R)]\subseteq C(R)[\lambda]$,\\
(2) $[QDer(R)_\lambda QC(R)]\subseteq QC(R)[\lambda]$,\\
(3) $[QC(R)_\lambda QC(R)]\subseteq QDer(R)[\lambda]$.
\end{lemma}

{\bf Proof.} Straightforward. \hfill $\square$

\begin{theorem}
Let $(R,\a)$ be a  multiplicative Hom-Lie conformal superalgebra. Then
\begin{eqnarray*}
GDer(R)=QDer(R)+QC(R).
\end{eqnarray*}
\end{theorem}

{\bf Proof.} For $D_{\mu}\in GDer_{\a^k}(R)$, there exist $D'_{\mu},D''_{\mu}\in \Omega$ such that
\begin{eqnarray}
[(D_{\mu}(a)_{\lambda+\mu}\alpha^{k}(b)]+(-1)^{|D||a|}[\alpha^{k}(a)_{\lambda}D'_{\mu}(b)]=D''_{\mu}([a_{\lambda}b]), \forall a,b\in R.
\end{eqnarray}

By (2.2) and (5.11), we get
\begin{eqnarray}
(-1)^{|D||b|}[\alpha^{k}(b)_{-\partial-\lambda-\mu}D_{\mu}(a)]+[D'_{\mu}(b)_{-\partial-\lambda}\alpha^{k}(a)]=D''_{\mu}([b_{-\partial-\lambda}a]).
\end{eqnarray}

By (2.1) and setting $\lambda'=-\partial-\lambda-\mu$ in (5.12), we obtain
\begin{eqnarray}
(-1)^{|D||b|}[\alpha^{k}(b)_{\lambda'}D_{\mu}(a)]+[D'_{\mu}(b)_{\mu+\lambda'}\alpha^{k}(a)]=D''_{\mu}([b_{\lambda'}a]).
\end{eqnarray}

Then,  changing the place of $a,b$ and replacing $\lambda'$ by $\lambda$ in (5.13) give
\begin{eqnarray}
(-1)^{|D||a|}[\a^{k}(a)_\lambda D_{\mu}(b)]+[D'_{\mu}(a)_{\lambda+\mu}\a^{k}(b)]=D''_{\mu}([a_{\lambda}b])
\end{eqnarray}
By (5.11) and (5.14), we have
\begin{eqnarray*}
&&[\frac{D_{\mu}+D'_\mu}{2}(a)_{\lambda+\mu}\a^{k}(b)]+(-1)^{|D||a|}[\a^{k}(a)_\lambda\frac{D_{\mu}+D'_\mu}{2}(b)]=D''_{\mu}([a_{\lambda}b]),\\
&&[\frac{D_{\mu}-D'_\mu}{2}(a)_{\lambda+\mu}\a^{k}(b)]-(-1)^{|D||a|}[\a^{k}(a)_\lambda\frac{D_{\mu}+D'_\mu}{2}(b)]=0.
\end{eqnarray*}
It follows that  $\frac{D_{\mu}+D'_\mu}{2}\in QDer_{\a^{k}}(R)$ and $\frac{D_{\mu}-D'_\mu}{2}\in QC_{\a^{k}}(R)$. Thus
\begin{eqnarray*}
D_\mu=\frac{D_{\mu}+D'_\mu}{2}+\frac{D_{\mu}-D'_\mu}{2}\in QDer_{\a^{k}}(R)+QC_{\a^{k}}(R),
\end{eqnarray*}
and we prove that $GDer(R)\subseteq QDer(R)+QC(R)$. The reverse inclusion relation follows from
(5.6) and Lemma 5.3.  \hfill $\square$

\noindent{\bf Theorem 5.5.} Let $(R,\a)$ be a  multiplicative Hom-Lie conformal superalgebra, $\a$   surjections
and $Z(R)$ the center of $R$. Then $[C(R)_{\lambda}QC(R)]\subseteq Chom(R,Z(R))[\lambda]$. Moreover, if $Z(R)=0$, then
$[C(R)_{\lambda}QC(R)]=0$.

{\bf Proof.} Since $\a$ are surjections, for any $b'\in R$, there exists $b\in R$ such that $b'=\a^{k+s}(b)$.
For any $D_{1\mu}\in C_{\a^{k}}(R), D_{2\mu}\in QC_{\a^{s}}(R), a\in R$, by (5.3) and (5.4), we have
\begin{eqnarray*}
&&[([D_{1\mu}D_2]_{\theta}(a))_{\lambda+\theta}b']\\
&=&[([D_{1\mu}D_2]_{\theta}(a))_{\lambda+\theta}\a^{k+s}(b)]\\
&=& [(D_{1\mu}(D_{2\theta-\mu}(a)))_{\lambda+\theta}\a^{k+s}(b)]-(-1)^{|D_1||D_2|}[(D_{2\theta-\mu}(D_{1\mu}(a)))_{\lambda+\theta}\a^{k+s}(b)]\\
&=& D_{1\mu}([D_{2\theta-\mu}(a)_{\lambda+\theta-\mu}\a^{s}(b)])-(-1)^{|D_1||D_2|+|D_2|(|D_1|+|a|)}[\a^{s}(b)(D_{1\mu}(a))_{\lambda+\mu}D_{2\theta-\mu}(\a^{k}(b))]\\
&=& D_{1\mu}([D_{2\theta-\mu}(a)_{\lambda+\theta-\mu}\a^{s}(b)])-(-1)^{|D_2||a|}D_{1\mu}([\a^{s}(a)_{\lambda}D_{2\theta-\mu}(b)]\\
&=&D_{1\mu}([D_{2\theta-\mu}(a)_{\lambda+\theta-\mu}\a^{s}(b)])-(-1)^{|D_2||a|}([\a^{s}(a)_{\lambda}D_{2\theta-\mu}(b)]\\
&=&0.
\end{eqnarray*}
Hence $[D_{1\mu}D_2](a)\in Z(R)[\mu]$, and then $[D_{1\mu}D_2] \in  Chom(R,Z(R))[\mu]$. If $Z(R)=0$, then $[D_{1\mu}D_2](a)=0$. Thus $[C(R)_{\lambda}QC(R)]=0$.
\hfill $\square$

\begin{proposition}
  Let $(R,\a)$ be a  multiplicative Hom-Lie conformal superalgebra, $\a$  surjections.  If $Z(R)=0$, then $QC(R)$ is a Hom-Lie conformal superalgebra if and only if $[QC(R)_{\lambda}QC(R)]=0$.
\end{proposition}

{\bf Proof. } $\Rightarrow$ Assume that $QC(R)$ is a Hom-Lie conformal superalgebra. Since $\a$ are surjections, for any $b'\in R$, there exists $b\in R$ such that $b'=\a^{k+s}(b)$.
For any $D_{1\mu}\in QC_{\a^{k}}(R), D_{2\mu}\in QC_{\a^{s}}(R), a\in R$, by (5.4), we have
\begin{eqnarray}
&&[\a^{k+s}([D_{1\mu}D_2]_{\theta}(a))_{\lambda+\theta}b']\nonumber\\
&=&[([D_{1\mu}D_2]_{\theta}(a))_{\lambda+\theta}\a^{k+s}(b)]+(-1)^{(|D_1|+|D_2|)|a|}[\a^{k+s}(a)_{\lambda}([D_{1\mu}D_2]_{\theta}(b))].
\end{eqnarray}
By (4.2) and (5.4), we obtain
\begin{eqnarray}
&&[([D_{1\mu}D_2]_{\theta}(a))_{\lambda+\theta}\a^{k+s}(b)]\nonumber\\
&=& [(D_{1\mu}(D_{2\theta-\mu}(a)))_{\lambda+\theta}\a^{k+s}(b)]-(-1)^{|D_1||D_2|}[(D_{2\theta-\mu}(D_{1\mu}(a)))_{\lambda+\theta}\a^{k+s}(b)]\nonumber\\
&=&(-1)^{|D_1|(|D_2|+|a|)}[\a^{k}(D_{2\theta-\mu}(a))_{\lambda+\theta-\mu}(D_{1\mu}(\a^s(b)))]\nonumber\\
&&-(-1)^{|D_1||D_2|+|D_2|(|D_1|+|a|)}[\a^{s}(D_{1\mu}(a))_{\lambda+\mu}(D_{2\theta-\mu}(\a^{k}(b)))]\nonumber\\
&=&(-1)^{|D_2||a|+|D_1|(|D_2|+|a|)}[\a^{k+s}(a)_{\lambda}(D_{2\theta-\mu}(D_{1\mu}(b)))]\nonumber\\
&&-(-1)^{|D_1||D_2|+|D_2|(|a|+|D_1|)+|D_2||a|)}[\a^{s}(a)_{\lambda+\mu}((D_{1\mu}(D_{2\theta-\mu}(b)))]\nonumber\\
&=& -(-1)^{(|D_2|+|D_1|)|a|)}[\a^{k+s}(a)_{\lambda}([D_{1\mu}D_2]_{\theta}(b))].
\end{eqnarray}
By (5.15) and (5.16)
\begin{eqnarray*}
[([D_{1\mu}D_2]_{\theta}(a))_{\lambda+\theta}b']=[([D_{1\mu}D_2]_{\theta}(a))_{\lambda+\theta}\a^{k+s}(b)]=0,
\end{eqnarray*}
and thus $[D_{1\mu}D_2]_{\theta}(a)\in Z(R)[\mu]=0$, since $Z(R)=0$. Therefore, $[QC(R)_{\lambda}QC(R)]=0$.

$\Leftarrow$ Straightforward.
\hfill $\square$
\begin{center}
 {\bf ACKNOWLEDGEMENT}
 \end{center}

  The paper is partially supported by  the NSF of China (No. 11761017), the Fund of Science and Technology Department of Guizhou Province (No. [2019]1021), the Key University Science Research Project of Anhui Province (Nos. KJ2015A294, KJ2014A183 and KJ2016A545),
  the  NSF of Chuzhou University (Nos. 2014PY08 and 2015qd01).

\renewcommand{\refname}{REFERENCES}


\begin{thebibliography}{99}


\bibitem{Ammar2010} F. Ammar,   A. Makhlouf,   Hom-Lie superalgerbas and Hom-Lie asmissible superalgebras,
J. Algebra, 324(2010), 1513-1528.

\bibitem{Ammar2013} F.  Ammar,  A. Makhlouf,  N. Saadaoui,  Cohomology of Hom-Lie superalgebras and $q$-deformed Witt superalgebra,
Czech. Math. J., 63(2013), 721-761.


\bibitem{Kac98} V. G. Kac, Vertex Algebras for Beginners,  Univ. Lect. Ser.  vol. 10 (Amer. Math. Soc., Providence, RI, 1996), second edition 1998.


\bibitem{Xu2000}  X. Xu, Quadratic conformal superalgebras, J. Algebra, 231 (2000), 1-38.

\bibitem{Yuan14} L. Yuan,  Hom Gel'fand-Dorfman bialgebras and Hom-Lie conformal algebras, J. Math. Phys., 55
(2014), 043507 .

\bibitem{Yuan17} L. Yuan, S. Chen,  C. He,  Hom Gel'fand-Dorfman super-bialgebras and Hom-Lie conformal superalgebras, Acta Mathematica Sinic, English Series, 33(1),
(2017), 96-116.

 \bibitem{Zhao2017} J. Zhao,   L. Chen,  L. Yuan,   Deformations and generalized derivations of
Lie conformal superalgebras,  J. Math. Phys., 58
(2017),  111702.


 \bibitem{Zhao2016} J. Zhao, L. Yuan,  L. Chen,   Deformations and generalized derivations of
Hom-Lie conformal algebras,   Science China Mathematics,  61 (2018), 797-812.

\end{thebibliography}
\end{document}